\documentclass[a4paper,12pt,oneside]{article}

\usepackage{eucal}
\usepackage{color}
\usepackage[all]{xy}
\usepackage[centertags]{amsmath}
\usepackage{latexsym}
\usepackage{stmaryrd}
\usepackage{amsfonts}
\usepackage{amssymb}
\usepackage{amsthm}
\usepackage{fancyhdr}
\usepackage{t1enc}
\usepackage[utf8]{inputenc} 
\usepackage[english]{minitoc}
\usepackage{fancyhdr}
\usepackage{newlfont}
\usepackage[active]{srcltx}
\usepackage{graphicx}
\usepackage{t1enc}
\usepackage{tikz}
\usepackage{makeidx}
\usepackage{appendix}
\usepackage[english]{babel}			
\usepackage{float}

\theoremstyle{plain}
\newtheorem{theorem}{Theorem}[section]
\newtheorem{remark}[theorem]{Remark}
\newtheorem{lemma}[theorem]{Lemma}
\newtheorem{exercise}[theorem]{Exercise}
\newtheorem{proposition}[theorem]{Proposition}
\newtheorem{definition}[theorem]{Definition}
\newtheorem{corollary}[theorem]{Corollary}
\newtheorem{notation}[theorem]{Notation}
\newtheorem{example}[theorem]{Exemple}
\newtheorem{prop}{Propri\'et\'e}[section]

\newcommand{\Ric}{\mathrm{Ric}}
\newcommand{\ric}{\widetilde{Ric}}

\newcommand{\C}{\mathbb C}

\newcommand{\beq}{\begin{equation}}
\newcommand{\eeq}{\end{equation}}
\newcommand{\beqn}{\begin{eqnarray}}
\newcommand{\eeqn}{\end{eqnarray}}
\newcommand{\bpro}{\begin{proposition}}
\newcommand{\epro}{\end{proposition}}
\newcommand{\blem}{\begin{lemma}}
\newcommand{\elem}{\end{lemma}}
\newcommand{\bdfn}{\begin{definition}}
\newcommand{\edfn}{\end{definition}}
\newcommand{\bcor}{\begin{corollary}}
\newcommand{\ecor}{\end{corollary}}
\newcommand{\bthm}{\begin{theorem}}
\newcommand{\ethm}{\end{theorem}}
\newcommand{\bex}{\begin{example}}
\newcommand{\eex}{\end{example}}
\newcommand{\brmq}{\begin{remark}}
\newcommand{\ermq}{\end{remark}}
\newcommand{\benum}{\begin{enumerate}}
\newcommand{\eenum}{\end{enumerate}}
\newcommand{\bitem}{\begin{itemize}}
\newcommand{\eitem}{\end{itemize}}
\newcommand{\bexer}{\begin{exercise}}
\newcommand{\eexer}{\end{exercise}}
\newcommand{\bproof}{\begin{proof}}
\newcommand{\eproof}{\end{proof}}
\newcommand{\eprop} {\end{prop} }
\newcommand{\bprop}{\begin{prop}}
\theoremstyle{plain}

\usetikzlibrary{arrows}
\usetikzlibrary{decorations.markings}
\tikzstyle directed=[postaction={decorate,decoration={markings,
    mark=at position .65 with {\arrow{latex}}}}]
\providecommand{\keywords}[1]
{
  \small	
  \textbf{Keywords:} #1
}
\providecommand{\thanks}[1]

\usepackage[english]{babel}
\title{Geometry of stable ruled surface over an elliptic curve}
\author{Arame Diaw\thanks{ supported by ANR-16-CE40-0008 project "Foliage"}\footnote{Univ Rennes, CNRS, IRMAR – UMR 6625, 35000 Rennes, France }}

\begin{document}
\maketitle						
\begin{abstract}
We consider the stable ruled surface $S_1$ over an elliptic curve.
There is a unique foliation on $S_1$ transverse to the fibration. The minimal self-intersection sections also define a $2$-web. We prove that the $4$-web defined by the fibration, the foliation and the $2$-web is locally parallelizable.\\

\end{abstract}
\bigskip
\tableofcontents
\section{Introduction}
Let $C$ be an elliptic curve on $\C$. In 1955, Atiyah proved in \cite{M.A.F} that, up to isomorphism, there are only two indecomposable ruled surfaces over $C$: the  semi-stable ruled surface $S_0\mapsto C$ and the stable ruled surface $S_1\mapsto C$. In this article, we study the geometry of the stable ruled surface. In fact, the surface $S_1$ can be seen as the suspension over $C$ of the unique representation onto the dihedral group $<-z,\dfrac{1}{z}>$ (see \cite{F.L}, page 23). Thus, we have a Riccati foliation $Ric$ on $S_1$ such that the generic leaf is a cover of degree $4$ over $C$ and it is the unique foliation transverse to the fibration. On the other hand, the holomorphic section $\sigma\colon C \mapsto S_1$ have self-intersection $\sigma.\sigma\geqslant 1$ and those having exactly $\sigma.\sigma=1$ form a singular holomorphic $2$-web $\mathcal{W}$. Finally, taking into account the fibration,
 we have a singular holomorphic $4$-web on $S_1$. The aim of this article is to study the geometry of this $4$-web composed by the Riccati foliation, the $2$-web $\mathcal{W}$ and the $\mathbb{P}^1$-fibration $\pi\colon S_1\mapsto C$.\\
Our first result is the following :
\bpro
 The discriminant $\Delta$ of the $2$-web $\mathcal{W}$ defined by the  $+1$ self-intersection sections on $S_1$ is a leaf of the foliation $\mathrm{Ric}$.
 \epro
 Using the isomorphism between the curve $C$ and its jacobian, we have the main result : 
 
 \begin{theorem} 
There exists a double cover $\varphi\colon C\times C\mapsto S_1$ ramified on $\Delta$ on which the lifted $4$-web $\textbf{W}$ is parallelizable.
\end{theorem}
 This $4$-web is locally comprised of pairwise parallel straight lines: its curvature is zero.\\
 Firstly, we show these results using only the properties of an elliptic curve and its jacobian and after, we use the theory of birational geometry to illustrate our results with computations on a trivialization  $S_1 \dashrightarrow C\times\mathbb{P}^1$.\\

This paper is part of my thesis work under the direction of Frank Loray and Frédéric Touzet.

\section{Preliminaries}
\subsection{Some properties on an elliptic curve}
 Let $C=\left\lbrace \left( x,y\right) \in{\mathbb{C}}^2,y^2=x\left( x-1\right)\left( x-t\right)  \right\rbrace \cup \{p_{\infty}\},$ 
 where $t\in\mathbb{C}\setminus \{0,1\}$ be an elliptic curve.
 Throughout this article, we use the following background of an elliptic curve. 
 \bpro\label{loi}
The set of points of $C$ forms an abelian group, with $p_{\infty}$ as the $0$ element and with addition characterized for any couple of points $p=\left(x_1,y_1\right)$, $q=\left(x_2,y_2\right)$ in $C$ by: 
\begin{enumerate}
 \item $-p=\left(x_1,-y_1\right)$ ;
 \item if $p\neq q,-q,$ then $p+q= \left(x_3,y_3\right)$ where $x_3=\lambda^2+(1+t)-x_1-x_2$, $y_3=\lambda(x_1-x_3)-y_1$ 
and $\lambda= \dfrac{y_2-y_1}{x_2-x_1}$ ;
 \item  if $p=q$ and $y_1\neq 0,$ then $2p= \left(\tilde{x},\tilde{y}\right)$ where 
 $\tilde{x}=\lambda^2-(1+t)-2x_1$ , $\tilde{y}= \lambda x_1-\tilde{x}-y_1$ and $\lambda=\dfrac{3x_1^2-2(1+t)x_1+t}{2y_1}$ ; 
 \item if $p=q$ and $y_1=0$ then $2p=p_{\infty}$. 
\end{enumerate}
  \epro
 \brmq
 The points $p_i=\left( i,0\right)$, where $i=0,1,t$ are the points of order $2$ on $C$ and the map
$$ \left\{
\begin{array}{cccc}
I\colon & C &\mapsto & C\\
 &\left( x,y\right)&\mapsto&\left( x,-y\right)
 \end{array} \right.$$
 is an automorphism of $C$ which fixes the points of order $2$: it is the standard involution of the curve $C$.
 \ermq
 If we denote $Jac(C)$, the jacobian of $C$, we have: 
 \blem
 There exists a bijection between $C$ and its jacobian defined by this following map: 
$$\left\{ \begin{array}{cccc} 
 &C&\longmapsto&Jac\left( C\right)\\
 &p &\longmapsto&[p]-[p_{\infty}]
\end{array} \right.$$
\elem
From now on, we will use the isomorphism between the additive group structure $(C,p_{\infty})$ and the group structure on $C$ induced by its jacobian.
\subsection{Ruled surface over an elliptic curve}
 Let $C$ be a smooth curve on $\C$.
 \bdfn
A ruled surface over $C$ is a holomorphic map of two dimensional complex variety $S$ onto  $C$ $\pi\colon S \mapsto C$ which makes $S$ a $\mathbb{P}^{1}$-fibration over $C$.
 \edfn
 \bex
The fiber bundle associated to a vector bundle of rank $2$ over $C$ is a ruled surface. We denote it, $\mathbb{P}(E).$
 \eex 
Conversely, we have the following theorem proved by Tsen in \cite{R.H} : 
\bthm\label{SR}
 Let $\pi\colon S \mapsto C$ be a ruled surface over $C$:
 \begin{enumerate}
 \item there exists a vector bundle $E$ of rank $2$ over $C$ such that $S=\mathbb{P}(E)$;
 \item there exists a section, i.e a map $\sigma\colon C \to S$ such that $\pi\circ\sigma=id$; 
 \item  $\mathbb{P}(E)\cong\mathbb{P}(E^{\prime})$ if and only if there is a holomorphic line bundle $L$ over $C$ such that $E\cong E^{\prime}\otimes L$.
\end{enumerate}  
\bdfn 
A ruled surface $\mathbb{P}(E)$ is decomposable if it has two disjoint sections.
 \edfn
 \ethm
The following lemma whose proof is in (\cite{M.M}, page $16$) shows the relationship between the ruled surface $S=\mathbb{P}(E)$ and the vector bundle $E$.
 \blem\label{friee}
There exists a one-to-one correspondance between the line subbundles of $E$ and the sections of $S$. Futhermore, if $\sigma_L$ is the section related to the line subbundle $L$ then: $${\sigma_L}.{\sigma_L}=\deg E-2\deg L$$ where $\deg(E)$ is the degree of the determinant bundle of $E$.
 \elem
 \begin{notation}
We recall that the notation $ {\sigma_L}.{\sigma_L}$ means the self-intersection of the section $\sigma_L$.
 \end{notation}
 \brmq
By the lemma \ref{friee}, $\mathbb{P}(E)$ is decomposable if and only if $E$ is decomposable, i.e $E= L_1\oplus L_2$ for line subbundles $L_i\subset E$.
 \ermq
 Consider $\kappa= \min\left\{\sigma.\sigma,\,\sigma\colon C\mapsto S\,/\pi\circ\sigma=id\right\}$. This number only depends on the ruled surface $S=\mathbb{P}(E)$. Indeed, it does not change when we replace $E$ by $E\otimes L$ for a line bundle $L$ on $C$.
 \bdfn
The ruled surface $\mathbb{P}(E)$ is stable if $\kappa>0$.
 \edfn

\bdfn
A minimal section of $S$ is a section $\sigma\colon C\to S$ such that, the self-intersection is minimal. That is to say, $\sigma.\sigma=\kappa$.
 \edfn
Using lemma \ref{friee}, we notice that a minimal section corresponds to a line subbundle of $E$ with maximal degree. Thus, the invariant $\kappa$ can be written as:
$$\kappa=\max\left\{\deg( E)-2\deg(L),\,L\hookrightarrow E_1\right\}$$
Now, we are interested in indecomposable ruled surfaces over an elliptic curve. Let  $\mathcal{O}_{C}\left( p_{\infty}\right)$ be the line bundle related to the divisor  $\left[p_{\infty} \right].$ There are unique nontrivial extensions of invertible sheaves:
  $$\xymatrix{ 
0\ar[r] &\mathcal{O}_C \ar[r] & E_0 \ar[r] & \mathcal{O}_{C} \ar[r] & 0
} $$ and
  $$\xymatrix{ 
0\ar[r] &\mathcal{O}_C \ar[r] & E_1 \ar[r] & \mathcal{O}_{C}\left( p_{\infty}\right)  \ar[r] & 0
} $$ Recall the following Atiyah's theorem as proved in (\cite{M.A}, Th.$6$.$1$):
  \bthm\label{Atiyah}
Up to isomorphism, the unique indecomposable ruled surfaces over $C$ are 
$S_0=\mathbb{P}(E_0)$  and $S_1=\mathbb{P}(E_1)$. 
  \ethm
 \brmq
Equivalency, any indecomposable vector bundle $E$ of rank $2$ on $C$ takes the form $E=E_i\otimes L$, for $i=0,1$ and $L$ a line bundle.  
 \ermq
 As our aim in this paper is the study of the ruled surface $S_1$, we will show firstly some important properties of $E_1$.
\blem
The degree of the maximal line subbundles of $E_1$ is zero.
\elem
\bproof
Let $L$ be a subbundle of $E_1$ and consider the quotient $M:=E_1/L$. By the following exact sequence $\xymatrix{ 
0\ar[r] &L \ar[r] & E_1 \ar[r] & M  \ar[r] & 0
} $ and the fact that $E_1$ is indecomposable, we have $H^1(M^{-1}\otimes L)\neq 0$.
Thus, due to Serre's duality, we have $2\deg(L)\leqslant \deg(E_1)$ and then $\deg(L)\leqslant 0$ because $\deg(E_1)=1$. Since the trivial line bundle $\mathcal{O}_C$ is a line subbundle over $E_1$, we have the result.
\eproof
\brmq
By this lemma, we can deduce that the ruled surface $S_1$ is stable. More precisely, up to isomorphism, it is the unique stable ruled surface over an elliptic curve. 
\ermq
If we consider $\max_{E_1}=\left\lbrace L\hookrightarrow E_1,\deg L=0\right\rbrace$ the set of line subbundles of $E_1$ having a maximal degree, we have:
\blem\label{cleee}
The jacobian of $C$ sets the parameters of the set $\max_{E_1}$. More precisely, the map
$$\left\{  \begin{array}{cccc}
 M\colon&\max_{E_1}&\longrightarrow&Jac\left( C\right)\\
 & L&\longmapsto&[L]\\
\end{array} \right.$$ is a bijection.
\elem
To prove this lemma, we have to use a key lemma of Maruyama in (\cite{M.M}, page $8$): 
\blem\label{clee}
Let $E$ be a vector bundle of rank $2$ over a curve. If $L_1$ and $L_2$
are distinct maximal line subbundles of $E$ such that 
$L_1$ and $L_2$ are isomorphic, then $E=L_1\oplus L_1.$
\elem
Now, we can prove lemma \ref{cleee}: 
\bproof
\begin{itemize}
\item Let $L_1$ and $L_2$ be two elements in $\max_{E_1}$ such that $L_1\cong L_2$. We have two possibilities, either $L_1=L_2$ or they are both distinct. According to the lemma \ref{clee}, the last case cannot occur because $E_1$ is not decomposable. Thus, the map $M$ is injective.
\item Let $L\in Jac\left( C\right)$ be distinct from the trivial line bundle. If we apply the functor  $Hom\left( L,--\right)$  to the exact sequence $$\xymatrix{ 
0\ar[r] &\mathcal{O}_C \ar[r]^f & E_1 \ar[r]^g & \mathcal{O}_{C}\left( p_{\infty}\right)  \ar[r] & 0
} $$ and we use Riemann Roch's theorem, we obtain  $\dim Hom\left( L\,,E_1\right)=1$. There exists a non zero morphism  $\tau\colon L\mapsto E_1.$ Thus, if we denote $D$ the effective divisor of zeros of $\tau$, then $L \otimes \mathcal{O}_C\left( D\right)$ is a line subbundle of $E_1$. Since $\deg\left( L\right) =0$,  $D$ is a effective divisor of zero degree, that is to say $\mathcal{O}_C\left( D\right)=\mathcal{O}_C$. Hence, $L$ is a line subbundle of $E_1$.
\end{itemize}
\eproof
\brmq
The minimal sections of $S_1$ have self-intersection equal to $1$ and 
they are parametrised by the jacobian which is isomorphic to $C$.
Hence for every point $\epsilon\in C$, we denote  $\sigma_{\epsilon}$ the minimal section corresponding via lemma \ref{cleee}
to the subbundle isomorphic to $\mathcal{O}_C\left( [\epsilon]-[p_{\infty}]\right)$.
\ermq
\blem\label{double clé}
Let $\sigma_{\epsilon}$ and $\sigma_{{\epsilon}^{\prime}}$ be two minimal sections of $S_1$. 
If we consider $\mathrm{D}$ their intersection divisor, we have:$$\pi(\mathrm{D})=[-\epsilon-\epsilon^{\prime}]$$ where 
$\pi\colon S_1 \mapsto C$ is the projection map.

\elem
\bproof
Intuitively, the divisor $D$ is defined by the points above at which the line bundles related to 
$\sigma_{\epsilon}$ and $\sigma_{{\epsilon}^{\prime}}$ coincide. More precisely, $\pi\left(\mathrm{D}\right)$
 is a effective divisor equivalent to divisor
$det(E_1)\otimes\mathcal{O}_{C}([p_{\infty}]-[{{\epsilon}}])
\otimes\mathcal{O}_{C}([p_{\infty}]-[{{\epsilon}^{\prime}}])$ which itself is equivalent to divisor 
$\mathcal{O}_{C}([-\epsilon-\epsilon^{\prime}])$. Since the degree of $D$ is equal to $1$, we obtain the result.
\eproof
\brmq\label{min}
Let $Q$ be a point of $S_1$ belonging to the fiber $\pi^{-1}(p)$.
If the minimal section $\sigma_{\epsilon}$ passes through the point $Q$, then the unique other minimal section passing through the same point 
$Q$ is the section $\sigma_{-p-\epsilon}$. They might be the same for some $Q$.
\ermq
\begin{figure}[H]
\centering
\includegraphics[scale=5, width=250pt]{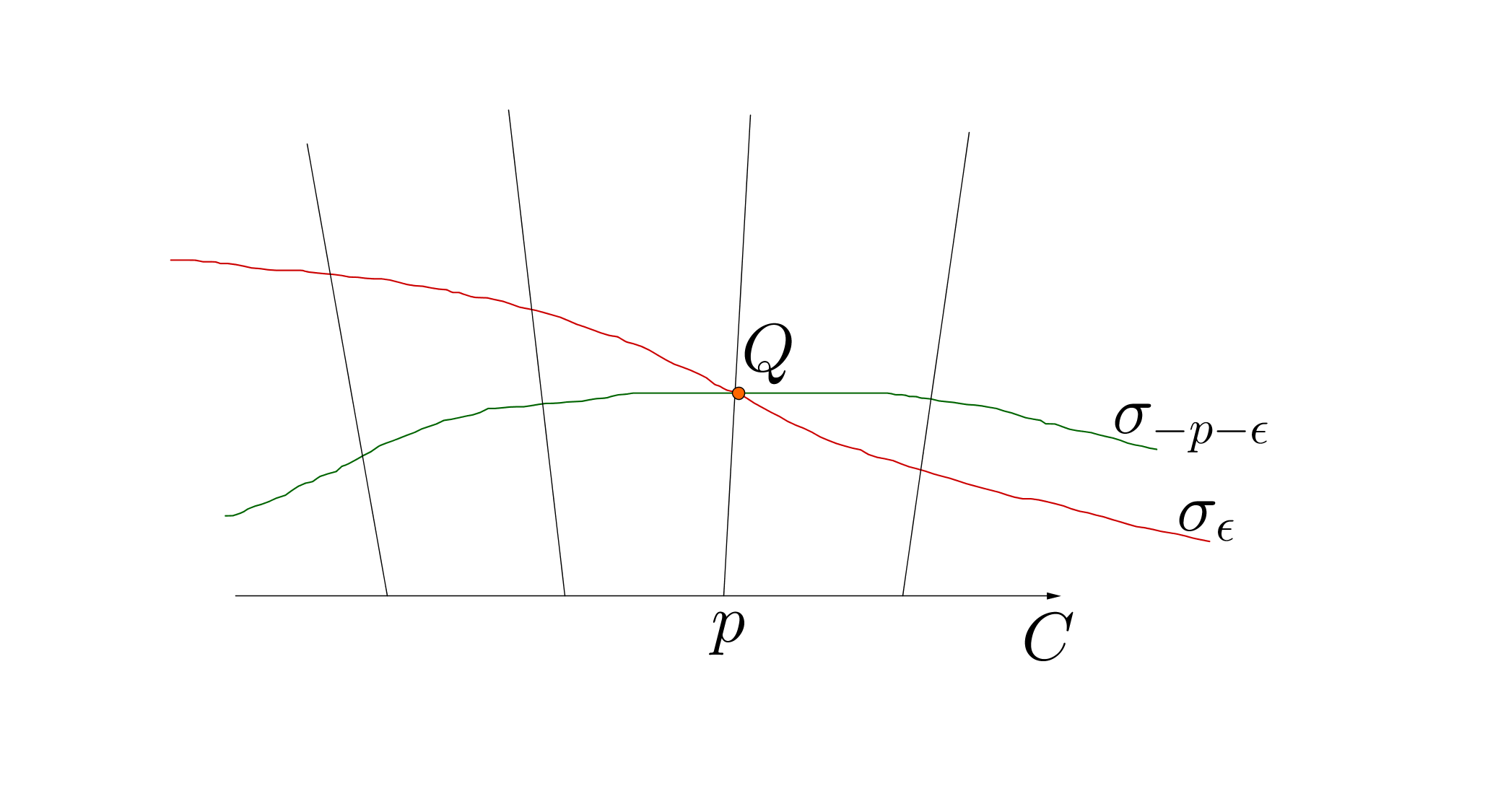}
\caption{intersection of two minimal sections on $S_1$}
\label{speci}
\end{figure}
We also have the following theorem proved by André Weil in \cite{A.W} :
 \bthm\label{André Weil}
 A holomorphic vector bundle on a compact Riemann surface is flat if and only if it is the direct sum of indecomposable vector
 bundles of degree $0$.
 \ethm
 By this theorem, the Atiyah's bundle $E_1$ is not flat because $\deg E_1=1.$ However what can we say about its associated ruled surface ?
 The answer of this question is given by Frank Loray and David Marin
in \cite{F.L}.\\
Consider $C$ as a torus $\mathbb{C}/\mathbb{Z}+\tau\mathbb{Z}$ and let   $\varrho\colon\mathbb{Z}+\tau\mathbb{Z}\to \mathbb{P}SL_2(\mathbb{C}) $ be the representation of the fundamental group of $C$ defined by $\varrho(1)=-z$ and $\varrho(\tau)=\dfrac{1}{z}$. Up to conjugacy in
$\mathbb{P}SL_2(\mathbb{C})$, it is the unique representation onto the $4$-order group $\Gamma=<-z,\dfrac{1}{z}>$.\\
The orbits of the elements of $\C\times\mathbb{P}^1$ modulo the action given by the representation and the universal cover of $C$ form a ruled surface over $C$, denoted by $\overline{E}$. It is obvious to see, the horizontal foliation of $\C\times\mathbb{P}^1$ lifts to a regular foliation $\mathcal{R}$ transverse to the fibration of $\overline{E}$. 
The foliated surface $(\overline{E},\mathcal{R})$ is called the suspension of $C$. Let $\mathcal{R}_{[x_0,z_0]}$ be a leaf passing through a the point $[x_0,z_0]$ of $\overline{E}$. Then, by definition, we have a isomorphism between $\mathcal{R}_{[x_0,z_0]}$ and the quotient $\C/G$ where $G=\{\alpha\in\mathbb{Z}+\tau\mathbb{Z} / \varrho_{\alpha}(z_0)=z_0\}$. Thus, we can deduce that every leaf of foliation $\mathcal{R}$ is a cover of $C$. It is not difficult to show that the intersection of any fiber of $\overline{E}$ and the leaf
$\mathcal{R}_{[x_0,z_0]}$ is given by the set $\{\varrho_{\alpha}(z_0), \alpha\in\mathbb{Z}+\tau\mathbb{Z}\}$.  Using the finitude of the representation, we have every leaf of the foliation is a cover of finite degree.  Moreover, the monodromy of the foliation $\mathcal{R}$ on any fiber is the representation $\varrho$.\\
Note that, the action of $<-z,\dfrac{1}{z}>$ in $\mathbb{P}^1$ gives two kind of orbits: orbits of order $4$ and three special orbits of order 2 given by $(-1,1)$, $(-i,i)$ and $(0,\infty)$. Therefore, the foliation $\mathcal{R}$ has a generic leaf which is a cover of order $4$ of $C$ and three special leaves, which is a cover of order $2$.
\bpro
The ruled surface $\overline{E}$ over $C$ is indecomposable such that its invariant $\kappa=1$. It is the ruled surface $S_1$. 
\epro
\bproof
If $\overline{E}$  is a decomposable ruled surface then its invariant $\kappa=0.$ Indeed, let $F$ be a leaf  of $\mathcal{R}$ and $\sigma_0$ a minimal section of $\overline{E}$, then we have: $F\equiv 4\sigma_0+bf$ or $F\equiv 2\sigma_0+b^{\prime}f$, where $f$ represents a fiber. Using the fact that $\overline{E}$ is decomposable, we can find a section $\sigma$ such that $\sigma.\sigma_0=0$. Since $F.\sigma\geqslant 0$ and $F^2=0$, we obtain that $\sigma_0.\sigma_0=0$.\\ Eventually, if we assume  that $\overline{E}$ is a decomposable ruled surface, we have $F\equiv 4\sigma_0$ or $F\equiv 2\sigma_0$ and then $F.\sigma=0$. The section $\sigma$ does not meet any leaf of $\mathcal{R}$, which does not make a sense because the foliation is regular.\\ The ruled surface $\overline{E}$ is then indecomposable. Hence, it is either isomorphic to $S_0$ or $S_1$.\\ By the same arguments above, $\overline{E}$ is not isomorphic to $S_0$, if not there would be a section which not intersects any leaf of $\mathcal{R}$.  
Thus, the ruled surface $\overline{E}$ is isomorphic to $S_1$ by Atiyah's theorem.\\
As up to conjugacy the representation $\varrho$ is the unique representation onto the $4$-order group $\Gamma=<-z,\dfrac{1}{z}>$, we deduce that the isomorphism between $\overline{E}$ and $S_1$ is the identity. 
\eproof
In summary, we have:
 \begin{theorem}
The ruled surface $S_1$ has a Riccati foliation $\mathrm {Ric}$ with  irreducible monodromy group $<-z,\dfrac{1}{z}>$.
  \end{theorem}
 \brmq
The foliation $\Ric$ has a generic leaf which is cover of degree $4$ over $C$ and three special leaves which are covers of degree $2$ over $C$. 
\ermq

\section{Geometry of the ruled surface  $S_1$ }
Let $\pi\colon S_1 \mapsto C$ be the stable ruled surface over $C$. 
\bpro\label{feuilles}
The automorphism group of $S_1$ is a group of order $4$ which is isomorphic to the $2$-torsion group in $C$.
\epro
\bproof
Let $\psi\colon S_1\mapsto S_1$ be a non trivial automorphism of $S_1$. Since the self-intersection is invariant by
automorphism $\psi$ preserves the set of $+1$ self-intersection sections on $S_1$. More precisely, for any
$\epsilon\in C$, there exists a unique point $r_{\epsilon}\in C$ such that
$\psi\left(\sigma_{\epsilon}\right)=\sigma_{r_{\epsilon}}$. The automorphism $\psi$ induces an automorphism
$\overline{\psi}$ of $C$ such that for any point $\epsilon\in C$ we have 
$\overline{\psi}\left(\epsilon\right)=r_{\epsilon}$. If we define $C$ as the complex torus 
$\mathbb{C}/\mathbb{Z}+\tau\mathbb{Z}$, we can write
for any $z\in \C,$ 
$\overline{\psi}\left(\overline{z}\right)=a\overline{z}+b$, where $a,b\in\mathbb{C}$ and
$a\left(\mathbb{Z}+\tau\mathbb{Z}\right)=\mathbb{Z}+\tau\mathbb{Z}$. \\ If we assume this automorphism has a fix point $\epsilon_0,$ then
by definition we have  $\psi\left(\sigma_{\epsilon_0}\right)=\sigma_{\epsilon_0}.$ Hence, using the lemma  
\ref{double clé}, we obtain that for any $p\in C,$ $\psi\left(\sigma_{-p-\epsilon_0}\right)=\sigma_{-p-{\epsilon_0}}$. For any
fiber, the automorphism $\psi$ is Moebius map which fixes at least three points: it is the trivial automorphism, which
does not make sense by hypothesis.\\ Therefore, the automorphism $\overline{\psi}$ has no fixed points, it is a translation like
$\overline{\psi}\left(\overline{z}\right)=\overline{z}+b$. As by definition we have: $\overline{\psi}\left(-p-\overline{z}\right)=-p-\psi(\overline{z})$, the point $b$ is a point of order $2$ of $C$.\\
Conversely, for any point $p_i$ of order $2$ on $C$,  we can define an automorphism $\Phi_{i}$, on
$S_1$  such that for any point $p\in C$, $\Phi_{i}$ restricts to the fiber $\pi^{-1}(p)$ is the unique Moebius map 
which associates the points of the sections $(\sigma_{p_\infty},\sigma_{p_0}, \sigma_{p_1},\sigma_{p_t})$ 
to the points of the sections $( \sigma_{p_\infty+p_i}, \sigma_{p_0+p_i},  \sigma_{p_1+p_i},  \sigma_{p_t+p_i})$ respectively. It is defined by:
$$\left\{  \begin{array}{cccc}
 \Phi_{i}\colon&S_1&\longrightarrow&S_1\\
 &z=\sigma_{\omega}\left( p\right) &\longmapsto&z^{\prime}=\sigma_{\omega+p_i }\left( p\right)\\
\end{array} \right.$$
There exists a one-to-one correspondance between the automorphisms of the fiber bundle $S_1$
and the points of order $2$ in $C$ which preserves the group structure. Hence we have : 
$$Aut_{C}\left(S_1\right)  = \left\lbrace \Phi_0, \Phi_1,\Phi_t,\Phi_{\infty}=Id\right\rbrace$$
\eproof
\bpro
The automorphism group of $S_1$ preserves the  foliation $\Ric$.
\epro
 \bproof
 Using the fact that the fundamental group of $C$ is abelian, we can extend the monodromy map over every fiber and
 regard it as automorphism on $S_1$ which fixes the basis $C$. Thus, we obtain that the monodromy group of the foliation $\Ric$
is a subgroup of order $4$ of $Aut_{C}\left(S_1\right)$ : they are isomorphic. The group $Aut_{C}\left(S_1\right)$ preserves the Riccati foliation on $S_1$.
 \eproof
 \bcor\label{monodromie}
The Riccati foliation $\Ric$ is the unique Riccati foliation on the ruled surface $S_1$.
 \ecor
 \bproof
Let $\mathcal{F}_1$ be a Riccati foliation on $S_1$. As its monodromy group is an abelian subgroup of $PGL(\C,2)$, we have three
possibilities for its monodromy representation :
 \begin{itemize}
\item  If the conjugacy class of the monodromy is the linear class defined by the group $\left\langle az\,,bz\right\rangle,$ there exists two disjoint invariant sections of $S_1$. Hence $S_1$ is a ruled surface related to the direct sum of two
line bundles over $C$. It does not make sense because $S_1$ is indecomposable.
\item If the conjugacy class of the monodromy is the euclidian class defined by the group $\left\langle z+1\,,z+s\right\rangle,$ there exists an invariant section on $S_1$ with zero self-intersection. In fact by the Camacho Sad's theorem (in \cite{C.S}), any invariant
curve of regular foliation has a zero self-intersection. This monodromy representation does not make sense in $S_1$ because
we have $\min\left\{\sigma.\sigma,\,\sigma\colon C\mapsto S_1\,/\pi\circ\sigma=id\right\}=1$.
\end{itemize}
The only remaining possibility is that the monodromy has image
the group  $<-z,\dfrac{1}{z}>$. Thus, the foliation
 $\mathcal{F}_1$ is conjugated to  $\Ric$ by an element in $Aut_{C}\left(S_1\right)$. As this automorphism group of the fibration $S_1$
 preserves the foliation $\Ric$, we have $\Ric=\mathcal{F}_1$. 
 
 \eproof
 \blem\label{+1}
There exists a ramified double cover of the ruled surface $S_1$ defined by the map : 
 $$ \left\{\begin{array}{cccc}
\varphi\colon & C\times Jac\left( C\right)&\longrightarrow & S_1\\
& \left(p,\epsilon\right)&\longmapsto& z=\sigma_{\epsilon}\left( p\right)
 \end{array}\right.$$ such that its involution is defined by :
 $$
 \left\{\begin{array}{cccc}
 \mathrm{i} \colon&C\times Jac\left( C\right)&\longrightarrow & C\times Jac\left( C\right)\\
& \left(p,\epsilon\right)&\longmapsto & \left( p\,,-p-\epsilon\right)
 \end{array}\right.$$        
 \elem
 \bproof
According to the lemma \ref{double clé}, three minimal sections cannot meet at the same point, then we deduce for any $p\in C$, the morphism
$$\left\{ \begin{array}{cccc}
 \varphi_p\colon&Jac\left( C\right) &\longrightarrow&\pi^{-1}\left( p\right)\\
 & \epsilon&\longmapsto&\sigma_{\epsilon}\left( p\right)\\
\end{array}\right.$$ is not constant: it is a ramified cover between Riemann surfaces.
Futhermore, by the remark \ref{min}, we know that at most two minimal sections can pass through a given point, then the map $\varphi$ is a ramified double cover.
 \eproof
 The immediate consequence of this lemma is the following :
\bthm\label{Tis}
There exists a singular holomorphic $2$-web $\mathcal{W}$ on $S_1$ 
defined by the minimal sections whose discriminant $\Delta$ is a leaf of the foliation $\Ric.$
\ethm
\bproof
 By lemma \ref{+1}, for any point $P\in \pi^{-1}\left( p\right)$ there exists a minimal section $\sigma_{r}$ passing through this point. Likewise, by lemma \ref{double clé} the minimal section 
 $\sigma_{-p-r}$ intersects transversally $\sigma_r$ at the point $P.$ As the sections  $\sigma_{r}$ and $\sigma_{-p-r}$ are
 distinct if and only if $2r\neq -p$, we deduce that there exists a singular holomorphic $2$-web on $S_1$ such that its discriminant is defined by : 
 $$\Delta=\cup_{p\in C}\left\lbrace P\in {\pi^{-1}}\left( p\right)\;/ P\in\sigma_{r},\;2r=-p\right\rbrace $$
In order to prove that $\Delta$ is a leaf of the Riccati, we need the following : 
\blem\label{clé}
There exists a linear foliation $\mathcal{F}$ on $C\times Jac\left( C\right)$ such that
$\varphi_{*}\mathcal{F}=Ric$.
\elem
\bproof
Assume that $C\times Jac\left( C\right)\simeq(\mathbb{C}/\mathbb{Z}+\tau\mathbb{Z})\times
(\mathbb{C}/\mathbb{Z}+\tau\mathbb{Z})$, and let be
 $(x,y)$ its local coordinates . If we consider the linear foliation
$\mathcal{\widetilde{F}}:=\mathrm{d} x+2\mathrm{d}y$ on $\C_{x}\times\C_{y}$,
then $\mathcal{\widetilde{F}}$ is invariant by the action of the lattice  $\mathbb{Z}+\tau\mathbb{Z}$. Thus we can lift
the foliation $\mathcal{\widetilde{F}}$ to a foliation, $\mathcal{F}$ on $ C\times Jac\left( C\right)$ 
such that the monodromy is defined by :
$$\left\{  \begin{array}{cccc}
 \xi\colon&\Lambda\longrightarrow&Aut(C)\\
 &\lambda\longrightarrow&z\mapsto z-\dfrac{1}{2}\lambda\\
\end{array} \right.$$
where $\Lambda$ is the lattice $\mathbb{C}/\mathbb{Z}+\tau\mathbb{Z}$.\\
The foliation $\mathcal{F}$ is transverse to the first projection on $C\times Jac\left( C\right)$ with a monodromy group isomorphic to the group of points of order $2$ $\left\lbrace p_{\infty},p_0,p_1,p_t\right\rbrace$.
Moreover, if $\mathcal{F}_{\left( p,\omega\right)}$ is the leaf passing through the point $\left( p,\omega\right)$,
then by definition we have : 
$$\mathrm{i}\left( \mathcal{F}_{\left( p,\omega\right)}\right) =\mathcal{F}_{\left( p, -p-\omega\right)}$$ 
where $i$ is the involution of the ramified double cover $\varphi$. Hence, $\varphi_{\ast}\mathcal{F}$ the direct image of the 
foliation $\mathcal{F}$ by $\varphi$ is a Riccati foliation on $S_1$ having the same monodromy group than $\Ric$. Using the 
uniqueness of $\Ric$ by the corollary \ref{monodromie}, we obtain the result.
\eproof
As by definition the curve $G=\left\lbrace\left( 2p,-p\right)/p\in C \right\rbrace $ is a leaf of the foliation $\mathcal{F}$,
using the foregoing lemma we can deduce that $\varphi\left( G\right) =\Delta$ is a leaf of $\Ric$. Which completes the proof of 
theorem \ref{Tis}.
\eproof
\subsection{Study of special leaves of the Riccati foliation $\Ric$}
According to lemma \ref{clé}, if $P=\sigma_{\omega}\left( p\right)\in S_1$ then the Riccati leaf passing through at this point is given 
by $$Ric_{P}=\left\lbrace\left( q,z\right) / z=\sigma_{\omega^{\prime}}\left( q\right) , 2\omega^{\prime}=2\omega+p-q\right\rbrace$$  
Thus, if we use this characterisation of the Riccati leaves on $S_1$, we have the following lemma :
 \blem\label{spéci}
 There exists three special leaves $Ric_0$, $Ric_1$ and $Ric_t$ of the foliation $Ric$ which are double cover of $C$. 
 More precisely, they are respectively the set of fixed points of the automorphisms $\Phi_0$, $\Phi_1$ and $\Phi_t$.  
  \elem
  \bproof
 We just give the proof for the leaf $Ric_0$ because it is the same process for the other special leaves.
  \begin{itemize}
 \item  Let $Ric_0$ be the Riccati leaf passing through the point $z_0=\sigma_{p_0}\left( p_0\right)$,
 then by definition, we have :
  $$Ric_0=\left\lbrace\left( p,z\right) / z=\sigma_{\omega}\left(  p\right) , 2\omega=p_{0}-p\right\rbrace $$
 According to the monodromy of $Ric$, if the leaf $Ric_0$ passes through the point 
 $z=\sigma_{\omega}\left(  p\right)$ then it passes through the points
  $\sigma_{\omega+p_i}\left(  p\right)$, where $p_i$ is a point of $C$ of order $2$. Since $2\omega=p_0-p$, we deduce from lemma \ref{double clé} that: 
 $\sigma_{\omega+p_0}\left(  p\right)=\sigma_{\omega}\left(  p\right)$ and 
 $\sigma_{\omega+p_1}\left(  p\right)=\sigma_{\omega+p_t}\left(  p\right)$ thus, $Ric_0$ meets any fiber of $S_1$ twice. It is a double cover over $C$.
\item Let $\Phi_0$ be the automorphism of $S_1$ related to the point $p_0$ defined by:
 $$\left\{  \begin{array}{cccc}
 \Phi_{0}\colon&S_1&\longmapsto&S_1\\
 &z=\sigma_{\omega}\left( p\right) &\longrightarrow&z^{\prime}=\sigma_{\omega+p_0 }\left( p\right)\\
\end{array} \right.$$ and consider the set of its fix points
$$\mathrm{Fix}_0=\left\lbrace \left( p,z\right) / \Phi_0\left( p,z\right) =\left( p,z\right) \right\rbrace$$ 
If $z=\sigma_{\omega}\left( p\right)$ is the fix point of $\Phi_0$, then by lemma \ref{double clé},
we have $2\omega=p_0-p$, and therefore $z\in Ric_0$. 
Conversely, if $z\in Ric_0$, we have $2\omega=p_0-p$, and 
according to lemma \ref{double clé}, we have, $\sigma_{\omega+p_0}\left( p\right)=  \sigma_{\omega}\left( p\right)$.
 Thus, we deduce that: 
$$Ric_0=\left\lbrace \left( p,z\right) / \Phi_0\left( p,z\right) =\left( p,z\right) \right\rbrace$$
\end{itemize}

\eproof
  According to all the foregoing, we have:
 \brmq 
The $2$-web given by the $+1$ self-intersection sections, the Riccati foliation and the $\mathbb{P}^1$-bundle
$\pi\colon S_1\mapsto C$ form a singular holomorphic $4$-web $\textbf{W}$ on $S_1$.
 \ermq
 \subsection{The geometry of the $4$-web $\textbf{W}$}
 Let $(x,y)$ be the local coordinates of $\mathbb{C}^{2}$. As the linear foliations $\mathcal{G}$ and $\mathcal{H}$ respectively defined by $\mathrm{d}y=0$ and by $\mathrm{d}y+\mathrm{d}x=0$ are
 invariant by the action of the lattice $\mathbb{Z}+\tau\mathbb{Z}$,
 we can lift them to a decomposable $2$-web $\mathcal{W^{\prime}}$ on $C\times Jac\left( C\right)$.
 \bpro
 The direct image $\varphi_{\ast}\left(\mathcal{W^{\prime}} \right)$ of the $2$-web $\mathcal{W^{\prime}}$ 
 by the ramified cover $\varphi$ is the $2$-web $\mathcal{W}$ on $S_1$ defined by the minimal sections.  
 \epro
 \bproof
 As by definition the $2$-web $\mathcal{W^{\prime}}$ is invariant by the involution of the ramified double cover $\varphi$, its direct image is also a singular holomorphic $2$-web on $S_1$.
Let $\left( p^{\prime},\,\epsilon^{\prime}\right)\in C\times Jac\left( C\right) $ and consider
\begin{itemize}
\item $A_{\epsilon^{\prime}}=\left\lbrace \left( p,\,\epsilon\right)\in 
 C\times Jac\left( C\right)/\,\epsilon=\epsilon^{\prime}\right\rbrace $
\item $B_{\epsilon^{\prime}}=\left\lbrace \left( p,\,\epsilon\right)\in 
 C\times Jac\left( C\right)/\,\epsilon=-p+\left( p^{\prime}+\epsilon^{\prime}\right) \right\rbrace ,$ 
 \end{itemize}
 the leaves of the $2$-web $\mathcal{W^{\prime}}$ passing through this point. 
 Since, using lemma \ref{double clé}, we have:
  $\varphi(A_{\epsilon^{\prime}})=\sigma	_{\epsilon^{\prime}}$ and  $\varphi(B_{\epsilon^{\prime}})=
 \sigma_{ -p^{\prime}-\epsilon^{\prime}}$, then the leaves of $\varphi_{\ast}\left(\mathcal{W^{\prime}} \right)$ are
 the minimal sections of $S_1$ which are the same along the discriminant $\Delta$.
 \eproof
The local study of the $4$-web $\textbf{W}$ on $S_1$ is the same as the
 $4$-web on $C\times Jac\left( C\right)$ given by the $2$-web $\mathcal{W^{\prime}}$,
 the foliation $\mathcal{F}$ and the $Jac\left( C\right)$-bundle defined by the first projection on 
 $C\times Jac\left( C\right)$. 
 \bthm\label{central}
Outside the discriminant locus $\Delta$, the $4$-web $\textbf{W}$ is locally  parallelizable .
 \ethm
 \bproof
According to the foregoing, the pull-back of $4$-web $\textbf{W}$ by the ramified cover $\varphi$ is locally
the $4$-web defined by $W\left( x, y, y+x, y+2x\right)$ on $\mathbb{C}^2$.  It is a holomorphic parallelizable web.
 \eproof
 \brmq
An immediate consequence of theorem \ref{central} is that the curvature of the $4$-web $\textbf{W}$ is zero.
 \ermq
The second part of this paper aims to use the theory of birationnal geometry in order to find the theoretic results of
the first part by computations on the birational trivialisation $C\times\mathbb{P}^{1}$.
 \section{Geometry of $4$-web $\textbf{W}$ after elementary transformations}
\begin{figure}[H]
\hspace*{-1.5cm}
\centering
\includegraphics[scale=3.5,width=500pt]{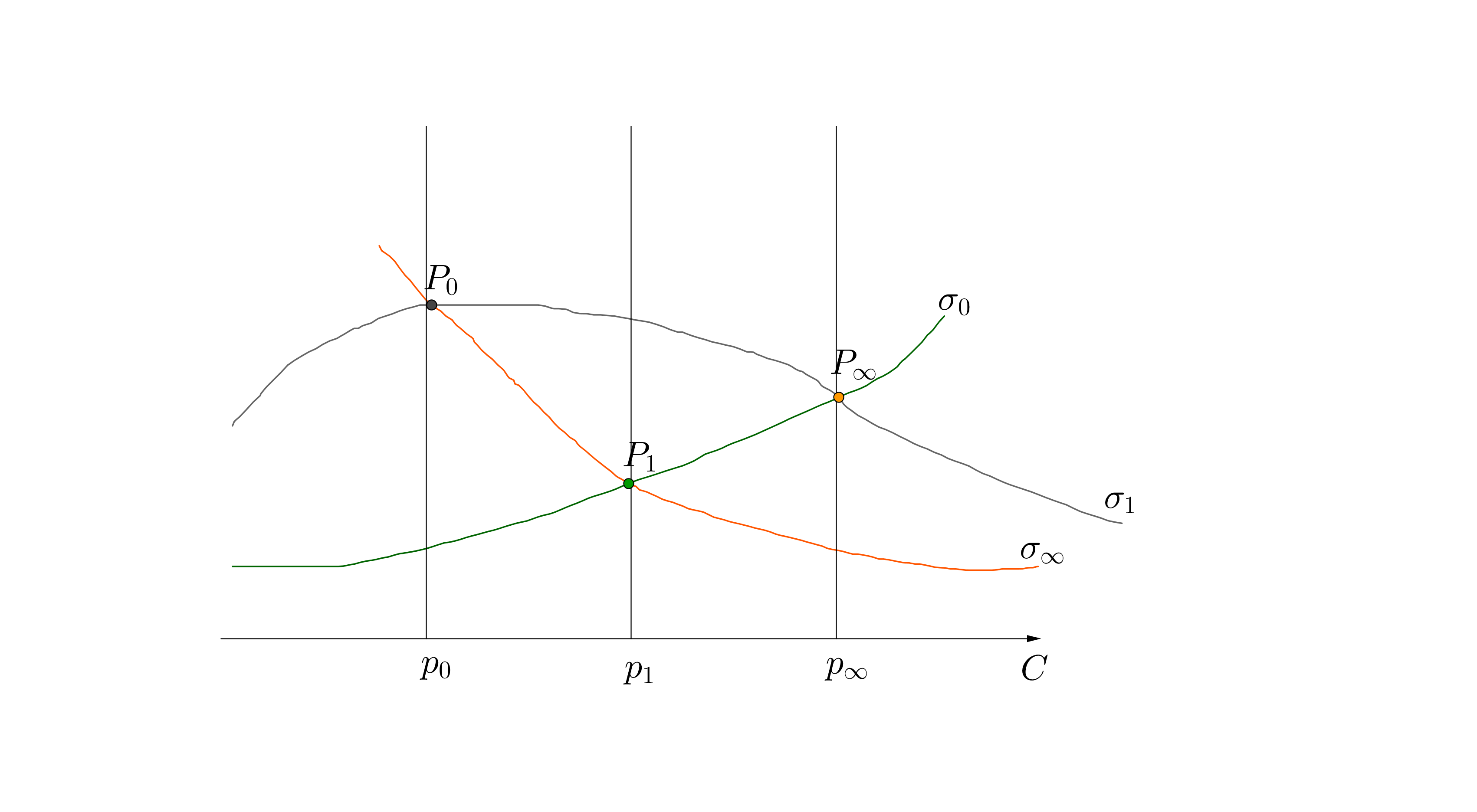}
\caption{Special points and special $+1$ self-intersection sections of the ruled surface $S_1$}
\label{s1}
\end{figure}

Let $\pi\colon S_1\mapsto C$ be the $\mathbb{P}^1$-bundle and 
$\left\lbrace p_0,p_1,p_t,p_{\infty}\right\rbrace$, the set of points of order $2$ in $C$. 
 \bdfn
An elementary transformation at the point $P\in\pi^{-1}(p) $ is the birational map given by the
composition of the blow-up of the point $P$, followed by the contraction of the proper transform of the fiber $\pi^{-1}(p)$.
 \edfn
\brmq
After elementary transformation at the point $P$, we obtain a new ruled surface with a point $\widetilde{P}$
which is the contraction of the proper transform of the fiber $\pi^{-1}(p)$.
 \ermq
How many elementary transformations do you need to trivialize the ruled surface $S_1$ ?
\blem 
The ruled surface $S_1$ is obtained after three elementary transformations 
 at the points $\widetilde{P}_0=\left(p_0,0\right) ,$
  $\widetilde{P}_1=\left( p_1,1\right) $ and  $\widetilde{P}_{\infty}=\left( p_{\infty},\infty\right) $ on the trivial bundle $C\times\mathbb{P}^1$.
  \elem
\begin{figure}[H]
\centering
\includegraphics[scale=3.5, width=500pt]{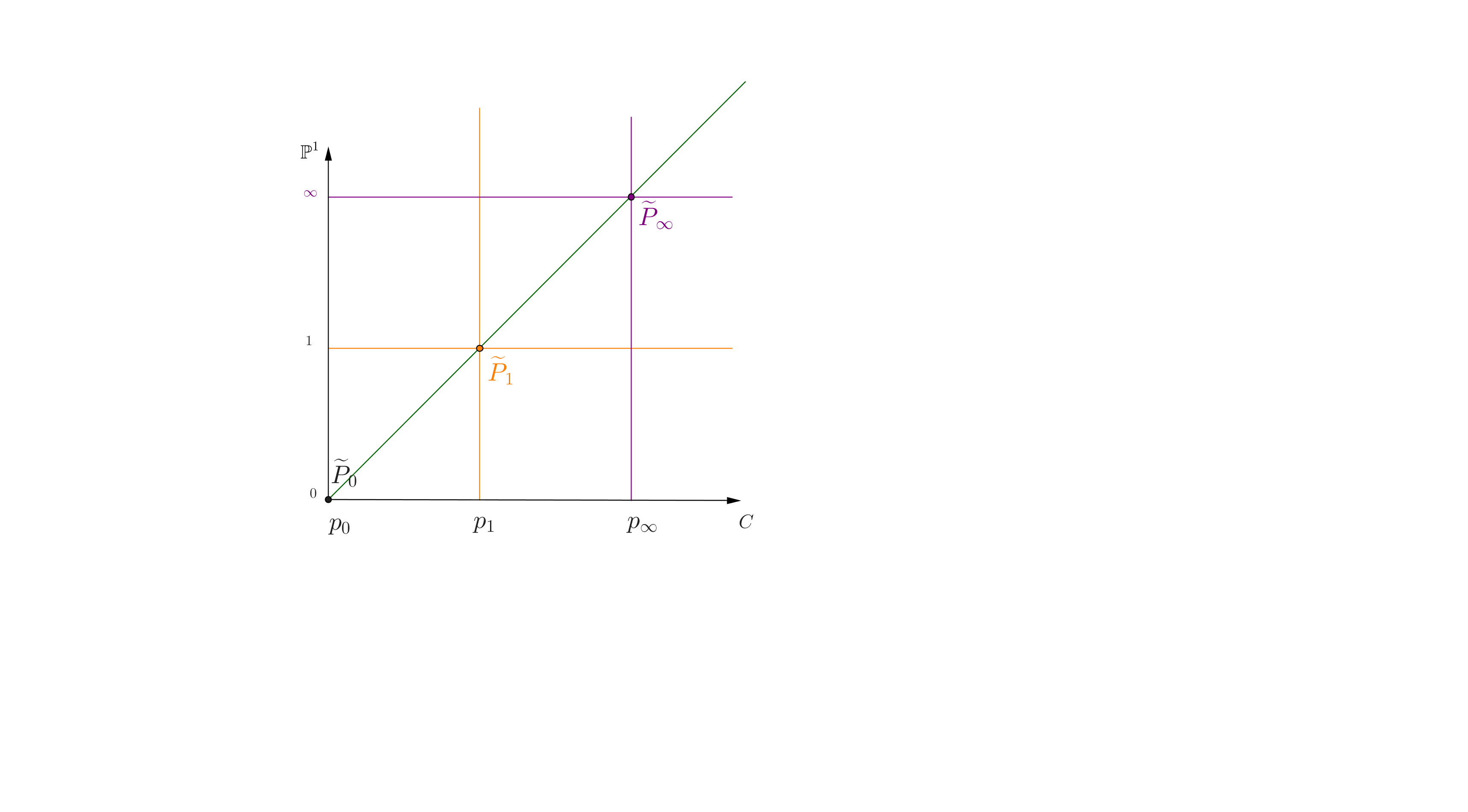}
\caption{Special points of the trivial $\mathbb{P}^1$-bundle over $C$}
\label{speci}
\end{figure}
 In fact, if we perform the elementary transformations of the three special points $\widetilde{P}_0$, $\widetilde{P}_1$ and $\widetilde{P}_{\infty}$ of $C\times\mathbb{P}^1$ : (see figure \ref{speci}),
we have a ruled surface $S$ with three special points $P_0$, $P_1$ and $P_{\infty}$ (see figure \ref{s1}).
Recall that after $3$ elementary transformations, if $\sigma$ is a section on $S$ such that $\sigma^{\prime}$ is its strict transform on the trivial bundle, we have:  $\sigma.\sigma= \sigma^{\prime}.\sigma^{\prime}+r$ where $r=\epsilon_0+\epsilon_1+\epsilon_{\infty}$ such that \begin{equation}
\epsilon_i=
\left\lbrace
\begin{array}{ccc}
-1  & \mbox{if} & P_i\in\sigma\\
+1 & \mbox{if} & P_i\notin\sigma \\
\end{array}\right.
\end{equation}
in particular, $r\in\{-1,1,-3,3\}$. Then, we can deduce that the ruled surface $S$ has a invariant $\kappa\leqslant 1$. 
\begin{enumerate}
\item If $\kappa=0$, then there is $\sigma.\sigma=0$ on $S$, then  its strict transform $\sigma^{\prime}.\sigma^{\prime}=odd$ on $C\times\mathbb{P}^1$. This cannot hold because all the sections of the trivial bundle have even self-intersection;
\item If $\kappa<0$, then there exists a $+2$ self-intersection section of $C\times\mathbb{P}^1$ passing through by three points. It is absurd because there exists a effective divisor equivalent to its normal bundle which contains at least three points. 
\end{enumerate}
According to these two cases, after our elementary transformations on the trivial bundle $C\times\mathbb{P}^1$, we obtain a ruled surface such that its invariant $\kappa=1$. Therefore, it is a stable ruled surface over an elliptic curve. \\
\subsection{The Riccati foliation on $S_1$ after elementary transformations}
 \bpro
After elementary transformations of the three special points $P_0$, $P_1$ and $P_{\infty}$ on $S_1$, the Riccati foliation $\Ric$ induces a Riccati foliation $\widetilde{Ric}$ on the trivial bundle $C\times\mathbb{P}^1$
such that the points $\widetilde{P}_0=\left(p_0,0\right) $, $\widetilde{P}_1=\left( p_1,1\right) $ and 
 $\widetilde{P}_{\infty}=\left( p_{\infty},\infty\right) $ are radial singularities. 
 \epro
 \bproof
 As the problem is local, we can prove it on the surface $\mathbb{C}^2$.  In our context, after elementary transformation at the origin,
a regular Riccati foliation becomes the pull-back of this holomorphic foliation $\dfrac{dz}{dx}=az^2+bz+c$ where $a,b,c\in \C$ by the birational map $\mathbb{C}^2\dashrightarrow\mathbb{C}^2$ ; $\left( x,z\right) \mapsto \left( x,xz\right) $.
 Thus, we obtain that a Riccati foliation such that the linear part looks like $xdz-zdx=0$. Therefore the origin is a radial 
 singularity. Finally, we can say the foliation $\Ric$ on $S_1$ is after elementary transformations a Riccati foliation
on $C\times\mathbb{P}^1$ having three radial singularities at the points
$\widetilde{P}_0$, $\widetilde{P}_1$ and $\widetilde{P}_{\infty}$.
 \eproof
If $((x,y),z)$ are cordinates of the trivial bundle $C\times \mathbb{P}^{1}$, then the foliation
$\widetilde{Ric}$ is defined by 
 $dz=\left[ a\left( x,y\right) z^2+b\left( x,y\right) z+c\left( x,y\right)\right] \dfrac{dx}{2y}$ 
 where $a,b,c$ are the meromorphic functions with pole of order $1$ at the points $p_0$, $p_1$ and $p_{\infty}$,
 i.e, 
 $a,b,c\in H^{0}\left( \mathcal{O}_{C}\left( p_0+p_1+p_{\infty}\right) \right)\simeq \mathbb{C}<1,\dfrac{1}{y},\dfrac{x}{y} >$.
It means that $a=\dfrac{a_0+a_{1}x+a_{2}y}{y}$, where $a_i$ are constant.
If we write the same relation for the functions  $b$ and $c$, we obtain that the foliation $\ric$ is locally defined by the following
$1$-form:
$$ydz=\left[ \left( a_0+a_{1}x+a_{2}y\right)z^2+ \left( b_0+b_{1}x+b_{2}y\right)z+\left( c_0+c_{1}x+c_{2}y\right)\right]\dfrac{dx}{2y}$$ 
As the foliation is invariant by the involution $\left(x,y \right)\mapsto \left(x,-y \right)$ on $C$, the coefficients
 $a_2,b_2,$ and $c_2$ are zero. Futhermore, if we use the relation on an elliptic curve,
 $y^2=x\left( x-1\right) \left( x-t\right)$, and the fact that the points 
 $\left( 0,0,0\right), \left( 1,0,1\right) ,\left(p_{\infty},\infty \right)$
are the radial singularities, we have $\ric$ is defined by the $1$-form :
$$w:=dz+\left[\dfrac{-z^2}{4x(x-1)}-\dfrac{z}{2x}+\dfrac{1}{4(x-1)}\right]dx$$
 \bpro
If we fix a point $(x_0,y_0)$ of $C$, the monodromy group of the foliation $\ric$ along of the fiber $\pi^{-1}(x_0,y_0)$ is an abelian group isomorphic
to a group
given by these automorphisms:
$\widetilde{\Phi}_0\colon z \mapsto \dfrac{z-x_0}{ z-1}$, 
$\widetilde{\Phi}_1\colon z \mapsto \dfrac{x_0}{z}$, 
$\widetilde{\Phi}_t\colon z \mapsto \dfrac{x_0\left( z-1\right)}{ z-x_0}$, $\widetilde{\Phi}_{\infty}\colon z\mapsto z$.
 \epro
  \bproof
Let $\widetilde{\sigma}_{\infty}:=\left\lbrace z=\infty\right\rbrace ,$ 
$\widetilde{\sigma}_{0}:=\left\lbrace z=0\right\rbrace $, $\widetilde{\sigma}_{1}:=\left\lbrace z=1 \right\rbrace$ and $\widetilde{\sigma}_{d}:=\left\lbrace z=x \right\rbrace $ be the four special sections obtained after elementary transformations. By definition of the monodromy group of $Ric$, we can see that for the point $p_0$ of order $2$, the automorphism $\widetilde{\Phi}_0$  restricted to any fiber is the unique Moebius transformation which relates respectively the points of the sections 
$(\widetilde{\sigma}_{0},\widetilde{\sigma}_{1},\widetilde{\sigma}_{\infty},\widetilde{\sigma}_{d})$ to the points of the sections $(\widetilde{\sigma}_{d}, \widetilde{\sigma}_{\infty}, \widetilde{\sigma}_{1},\widetilde{\sigma}_{0})$  . Using the same process for the other automorphisms, we obtain the result.

   \eproof 
 We can also describe the special leaves of the foliation $\ric$.
In fact, if we consider
$\phi_{i}\colon C\times\mathbb{P}^1\mapsto C\times\mathbb{P}^1$; $(x,y,z)\mapsto (x,y,\widetilde{\Phi_i}(z))$, then according to 
lemma \ref{spéci}, the special leaves are defined by :
\begin{enumerate}
\item $\widetilde{Ric_0}:=\left\lbrace \left( x,y,z\right) , \phi_0\left( x,y,z\right)=(x,y,z)\right\rbrace =\left\lbrace\left( x,z\right) , -z^2-x+2z=0 \right\rbrace$ 
\item $\widetilde{Ric_1}:=\left\lbrace \left( x,y,z\right) , \phi_1\left( x,y,z\right)=(x,y,z)\right\rbrace =\left\lbrace\left( x,z\right),-z^2+x=0 \right\rbrace$ 
\item $\widetilde{Ric_t}:=\left\lbrace \left( x,y,z\right) , \phi_t\left( x,y,z\right)=(x,y,z)\right\rbrace =\left\lbrace\left( x,z\right) ,  z^2-2xz+x=0 \right\rbrace$ 
\end{enumerate}
Now it is natural to ask if we can find the expression of the leaf of order $4$ of $\ric$. To do this, we use the special
leaves to find a first integral.
Let $$f_0 := -z^2+2z-x, \quad f_1 := -z^2+x,\quad f_t := z^2-2xz+x$$ be the polynomials 
which define respectively the leaves $\widetilde{Ric_0}$, $\widetilde{Ric_1}$, $\widetilde{Ric_t}$ 
and consider the function $\gamma\colon C_x\times\mathbb{P}^1_{z}\mapsto C_x\times\mathbb{P}^1_{y}$  ;
$\left(x,z \right) \mapsto
 \left( x,F\left( x,z\right)\right)$, where 
$$F\left( x,z\right)=\dfrac{xf_{0}^2}{xf_{0}^2-(x-1)f_{1}^2}$$
The pull-back $\gamma^{\ast}dy$ of the $1$-form $dy$ by $\gamma$ is a foliation on $C_x\times\mathbb{P}^1_{z}$ having the function $F\left(x,z\right)$ as a first integral 
 and such that the curves  $\widetilde{Ric_0}$, $\widetilde{Ric_1}$ and $\widetilde{Ric_t}$ are invariant. Hence, it is the foliation
$\ric$.
 we deduce that :
 \blem
 The foliation $\ric$ on $C\times\mathbb{P}^1$ has a rational first integral defined by the following function :
  $$F\left( x,z\right)=\dfrac{x(z^2-2z-x)^2}{(-z^2+2xz-x)^2}$$
  \elem
\subsection{The generic $2$-web after elementary transformations}
\begin{figure}[H]
\hspace*{1.5cm}
\centering
\includegraphics[scale=3.5, width=360pt]{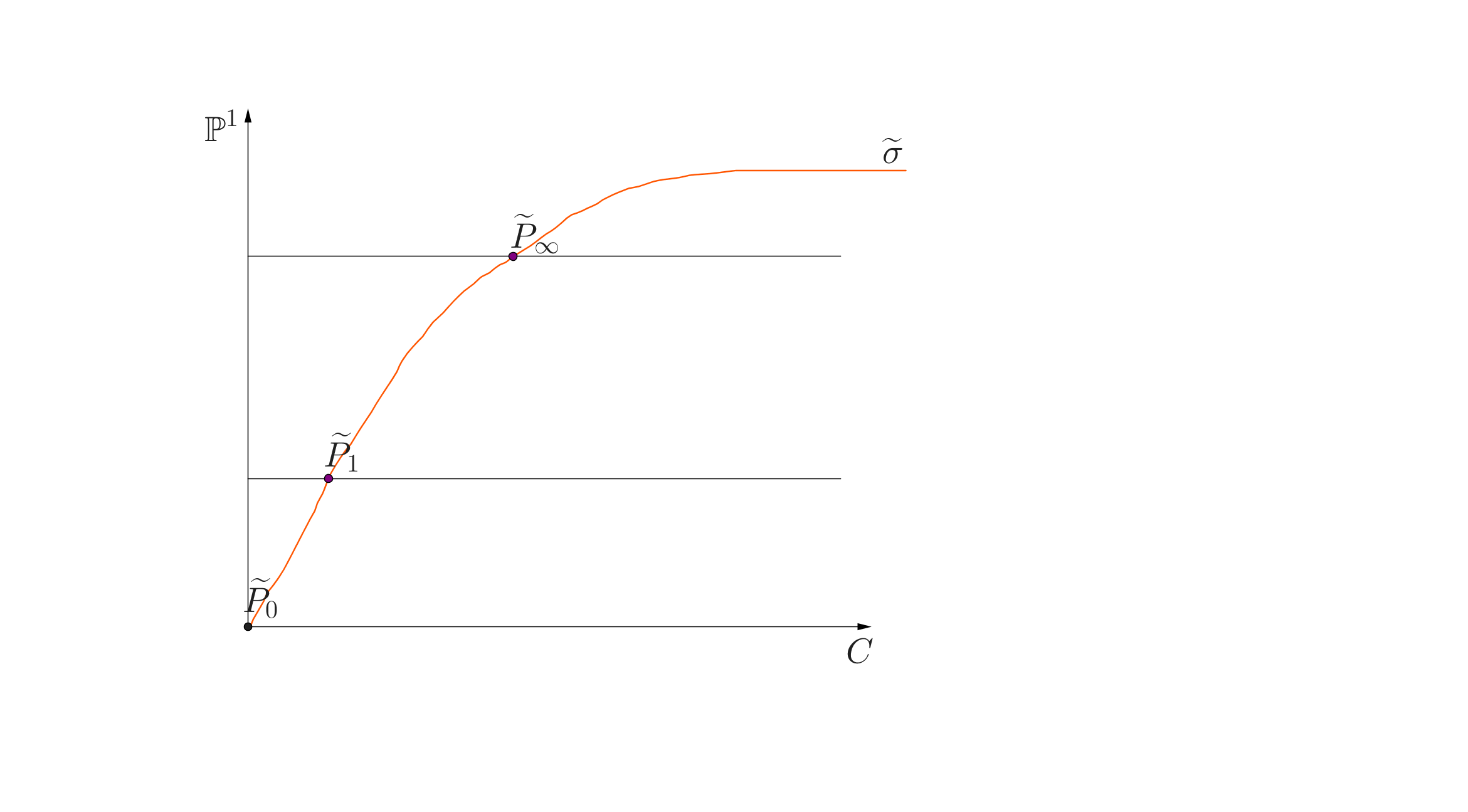}
\caption{Generic $+4$ self-intersection section}
\label{+4}
\end{figure}

After elementary transformations at the three special points on $S_1$, the generic $+1$ self-intersection sections (i.e not passing through the three special points) become the $+4$ self-intersection sections of $C\times\mathbb{P}^1$ passing through the points $\left( 0,0,0\right)$, $\left( 1,0,1\right) $ and $\left( p_{\infty},\infty\right)$: see figure \ref{+4}.\\
 \blem
A $+4$ self-intersection section passing through the points $\widetilde{P}_0$, $\widetilde{P}_1$ and $\widetilde{P}_{\infty}$ is either given by
the graph $z=\dfrac{(1-x_0)(y_{0}x-x_{0}y)}{y_0(x-x_0)}$, or the graph $z=x$.
    \elem    
\bproof
 If $\sigma\colon C\mapsto\mathbb{P}^1$ is a $+4$ self-intersection section on the trivial bundle, then it defines
 a rational map of degree $2$ generated by two sections $\sigma_1$ and $\sigma_2$ of a line bundle of degree $2$
 over $C$; more precisely, for any point $\left( x,y\right) \in C$,
 $\sigma\left( x,y\right) =\left( \sigma_1\left( x,y\right) : \sigma_2\left( x,y\right) \right)$.
 Since for any line bundle of degree $2$ over $C$, there exists a point $p=\left( x_{0},y_{0}\right) \in C$
 such that $L=\left[ p\right]+\left[ p_{\infty}\right]$, we have two cases: 
\begin{itemize}
\item if $p\neq p_{\infty}$, according to the Riemann Roch's theorem, 
$H^0\left( L\right) =\mathbb{C}\left\langle  y-y_{0},x-x_{0}\right\rangle $ and then,  $\sigma$ is a graph given by 
$$z=\dfrac{a\left( y-y_{0}\right) +b\left(  x-x_{0}\right) }{c\left( y-y_{0}\right) +d\left(  x-x_{0}\right)},\quad a,b,c,d\in C$$
Using the fact that the section passes through the points $\widetilde{P}_0$, $\widetilde{P}_1$ and $\widetilde{P}_{\infty}$ and the puiseux parametrisation of elliptic curve at the infinity point is given by $t\mapsto (\dfrac{1}{t^2},\dfrac{1}{t^3})$, we obtain 
a system of equations which solutions are $ \left\lbrace a = d\dfrac{x_0(x_0-1)}{y_0}, b = -d(x_0-1), c = 0, d = d\right\rbrace$
where $d\neq 0$;
 \item if $p=p_{\infty}$ then $H^0(L)=\mathbb{C}<1,\dfrac{1}{x}>$, likewise using the fact that the section
 passes through the points $\left( 0,0,0\right)$, $\left( 1,0,1\right) $ and $\left( p_{\infty},\infty\right)$, we obtain that
 $\sigma$ is the graph $z=x$. 
\end{itemize}
  \eproof
From now on, unless otherwise mentionned, we will assume that a $+4$ self-intersection section is a section which
passes through the points $\widetilde{P}_0$, $\widetilde{P}_1$ and $\widetilde{P}_{\infty}$.\\
 Now using the birational trivialisation of $S_1$, we can give another proof to show
that the minimal sections of $S_1$ define a singular $2$-web and its discriminant is a leaf
 of the Riccati foliation on $S_1$.
 \bpro
 For any point of  $C\times\mathbb{P}^1$, there exists a $+4$ self-intersection section which passes 
 through this point.
 \epro
 \bproof
 Let $\left( u,v,z\right)\in C\times\mathbb{P}^1$ such that $v\neq 0$, we have to find the points
 $\left( x_0,y_0\right)\neq(u,v)$ of $C$ such that  $z=\dfrac{(1-x_0)(y_{0}u-x_{0}v)}{y_0(u-x_0)}$.\\
 Using the fact that $y_0^2 = x_0(x_0-1)(x_0-t)$ and $ v^2 = u(u-1)(u-t)$, we have the following equation:\\
 $(\star): (u-z)^2x_0^3-[(2(uz-u))(-z+u)-(-z+u)^2t-v^2]x_0^2+[(uz-u)^2-(2(uz-u))(-z+u)t+v^2]x_0-(uz-u)^2t=0$
 \begin{enumerate}
 \item if $(u-z)=0$, then $(\star)$ becomes $(x_0-u)\left(x_0-\dfrac{t(u-1)}{u-t}\right)=0$.
 As by hypothesis $v\neq 0$, we obtain two solutions given by the point $(x_0,y_0)$ such that $x_0=\dfrac{t(u-1)}{u-t}$ and the point $p_{\infty}$; 
 \item if $(u-z)\neq 0$, then the solutions verify the following second degree equation:
 $$(\bigstar):(u-z)^2x_0^2+[(-t-u)z^2+2u(t+1)z-u(t+u)]x_0+tu(z-1)^2=0$$
 \end{enumerate}
 \eproof
 The $+4$ self-intersection sections define a singular holomorphic $2$-web $\mathcal{W}$ such that the discriminant $\Delta$ is the discriminant of the equation $\bigstar$. Thus we have :
$$\Delta:= (t-u)z^4-4(t-1))uz^3+2u(2tu+t-u-2)z^2-4u^2(t-1)z+u^2(t-u)=0$$
 \blem
 The discriminant of the $2$-web $\mathcal{W_1}$ is a leaf of order $4$ of the Riccati foliation $\ric$.
 \elem
 \bproof
In fact, by the definition of the first integral of the foliation $\ric$, we have: 
$$F(x,z)-t=-\dfrac{(t-u)z^4-4(t-1))uz^3+2u(2tu+t-u-2)z^2-4u^2(t-1)z+u^2(t-u)}{(2xz-z^2-x)^2}$$
Therefore, the first integral is constant along of the discriminant 
$\Delta$.
 \eproof
According to the foregoing, on the birational trivialisation of $S_1$, we have a $4$-web $\mathcal{W}_4$ defined by 
the $2$-web $\mathcal{W_1}$, the Riccati foliation $\ric$ and the trivial fiber bundle. 
\subsection{Geometry of the $4$-web $\mathcal{W}_4$}
We want to find the slopes of the leaves of $\mathcal{W}_4$ in order to represent it by a differential 
equation.
Let $\left( x_{0},y_{0},z_0\right)\in C\times\mathbb{P}^1$ be a generic point.
As the leaves of the $2$-web $\mathcal{W_1}$ passing through this point are respectively the graph $z=\dfrac{(1-a_0)(b_{0}x-a_{0}y)}{b_0(x-a_0)}$ and  $z=\dfrac{(1-{a_0}^{\prime})({b_{0}}^{\prime}x-{a_{0}}^{\prime}y)}{{b_0}^{\prime}(x-{a_0}^{\prime})}$
such that the points $\left( a_0, b_0\right)$ and $\left( {a_{0}}^{\prime},{b_{0}}^{\prime}\right) $ verify
the equation $(\bigstar)$, we deduce that their slopes at the point $\left( x_{0},y_{0},z_0\right)$ are respectively
given by the following formulas:
\begin{enumerate}
\item $Z_1=\dfrac{1-a_0-z_0}{x_0-a_0}+\left(z_0+\dfrac{(a_0-1)x_0}{x_0-a_0}\right)\left( \dfrac{ 3x_0^2-2\left( 1+t\right)x_0+t}{2x_0\left(x_0-1\right) \left( x_0-t\right)}\right)  $;
\item $Z_2= \dfrac{1-{a_0}^{\prime}-z_0}{x_0-{a_0}^{\prime}}+\left(z_0+\dfrac{({a_0}^{\prime}-1)x_0}{x_0-{a_0}^{\prime}}\right)\left( \dfrac{ 3x_0^2-2\left( 1+t\right)x_0+t}{2x_0\left(x_0-1\right) \left( x_0-t\right)}\right) $.
\end{enumerate}
Thus, the $2$-web $\mathcal{W_1}$ is defined by the following differential equation: 
$$\left(\dfrac{\mathrm{d}z}{\mathrm{d}x}\right)^2-\left(\dfrac{z^2+2(x-1)z-x}{2x(x-1)}\right)\dfrac{\mathrm{d}z}{\mathrm{d}x}+\dfrac{z(z-1)((2tx-x^2-t)z-x^3+x^2-tx+2)}{4x^2(x-1)^2(t-x)} $$
Futhermore, if we consider
$$Z_0=\dfrac{1}{4}\left[\dfrac{ \left( z_{0}^2+2(x_0-1)z_0-x_0\right)}{x_0\left( x_0-1\right)}\right],$$
the slope of the foliation $\ric$ at the point $\left( x_{0},y_{0},z_0\right)$, then the $4$-web 
$\mathcal{W}_4$ is locally equivalent to the $4$-web on the complex plane given by
$\mathcal{W}\left( \infty,Z_0,Z_1,Z_2\right)$.
 \bthm
The $4$-web $\mathcal{W}\left( \infty,Z_0,Z_1,Z_2\right)$ is locally equivalent to a parallelizable $4$-web. 
 \ethm
 \bproof
The pull-back of the foliation $\ric$ by the multiplication of order $2$ on $C$ is another Riccati foliation $\ric_2$ on $C\times \mathbb{P}^1$ with trivial monodromy. \\
Let $M_2\colon C\mapsto C,$ be the multiplication of order $2$ on $C$ then, for any point $\left( x,y\right)\in C$
the first projection of $M_2\left( x,y\right)$ is given by the following formula :
 $$pr_1\circ M_2\left( x,y\right)=\dfrac{ (3x^2-2(t+1)x+t)^2}{4x(x-1)(x-t)}+\left( 1+t\right) -2x$$
Using the pull-back of the special leaves, we can choose three curves by:
\begin{enumerate}
\item $C_0:= \left\lbrace \left(x,y,z\right) ,z=z_0=\dfrac{-x^2+t}{2(t-x)}\right\rbrace $
 \item $C_1:= \left\lbrace \left(x,y,z\right) ,z=z_1=\dfrac{(-x^2+t)}{2y}\right\rbrace $
 \item  $C_2:= \left\lbrace \left(x,y,z\right) ,z=z_2=-\dfrac{(-x^2+t)}{2y}\right\rbrace $
    \end{enumerate}
 which are the leaves of $\ric_2$. Now, if we consider the map  
 $\psi\colon C\times\mathbb{P}^1\mapsto C\times\mathbb{P}^1$ which for any coordinate $(X,Z)$ relates:
 $$\psi(X,Z)=\left(\dfrac{ (3X^2-2(t+1)X+t)^2}{4X(X-1)(X-t)}+\left( 1+t\right) -2X,
 \dfrac{Z\mu z_1-z_0}{Z\mu-1}\right)$$ where $\mu=\dfrac{z_3-z_0}{z_3-z_1}$, then the pull-back of the first
 integral of $\ric$ by $\psi$  is the following meromorphic function :
 $$(\psi_{\ast}F)\left( X,Z\right)=\dfrac{(Z^2-2Z+2)^2}{Z^2(Z-2)^2 }$$
Finally, the foliation $\psi_{\ast}\ric$ is locally defined by the $1$-form $dZ=0$. Likewise, the pull-back of the 
slopes $Z_1$ and $Z_2$ by $\psi$ defines a $2$-web such that the leaves verify the following differential equation :

$$(\star \star): \left( \dfrac{dZ}{dX}\right)^2+\dfrac{(t-1)Z^4+(-4t+4)Z^3+(4t-8)Z^2+8Z-4}{4X(X-1)(X-t)}=0$$
In summary, the $4$-web $\psi_{\ast}\mathcal{W}\left( \infty,Z_0,Z_1,Z_2\right)$ is locally equivalent to the web $\mathcal{W}\left( \infty,0,\beta,-\beta\right)$, where $\beta$ is a solution of  $(\star \star).$
As the $4$-web $\mathcal{W}\left( \infty,0,\beta,-\beta\right)$ has
a constant cross-ratios equal to $-1$ and all the $3$ subweb are hexagonal, we can deduce that it is locally parallelizable.
 \eproof
\nocite{ripoll}
 \nocite{R.F}
 \nocite{M.A}
 \nocite{M.S}
 \nocite{L.P}
\bibliographystyle{plain}

\bibliography{biblio}

\end{document}